\documentclass[fleqn]{mat01}
\usepackage{times,mathtimy,amssymb,latexsym}
\begin{document}

\setcounter{page}{267} \firstpage{267}

\newtheorem{theore}{Theorem}
\renewcommand\thetheore{\arabic{theore}}
\newtheorem{theor}[theore]{\bf Theorem}

\newtheorem{propo}{\rm PROPOSITION}
\newtheorem{coro}{\rm COROLLARY}

\def\scase{\trivlist \item[\hskip \labelsep{\it Special cases}]}
\def\efros{\trivlist \item[\hskip \labelsep{\bf Efros' theorem \cite{1}.}]}
\def\lem{\trivlist \item[\hskip \labelsep{\it Lemma.}]}

\renewcommand\theequation{\arabic{section}.\arabic{equation}}

\title{Fractional extensions of some boundary value problems\\ in oil
strata}

\markboth{Mridula Garg and Alka Rao}{Fractional extensions of some
boundary value problems in oil strata}

\author{MRIDULA GARG and ALKA RAO}

\address{Department of Mathematics, University of Rajasthan,
Jaipur~302~004, India\\
\noindent E-mail: gargmridula@gmail.com;
alkarao.ar@rediffmail.com}

\volume{117}

\mon{May}

\parts{2}

\pubyear{2007}

\Date{MS received 27 April 2006; revised 16 November 2006}

\begin{abstract}
In the present paper, we solve three boundary value problems
related to the temperature field in oil strata~-- the fractional
extensions of the incomplete lumped formulation and lumped
formulation in the linear case and the fractional generalization
of the incomplete lumped formulation in the radial case. By using
the Caputo differintegral operator and the Laplace transform, the
solutions are obtained in integral forms where the integrand is
expressed in terms of the convolution of some auxiliary functions
of Wright function type. A~generalization of the Laplace transform
convolution theorem, known as Efros' theorem is widely used.
\end{abstract}


\maketitle

\section{Introduction}

The heat equation for a porous medium is derived in \cite{1,10}
under several assumptions on the model. Besides the exact
formulation of the temperature field in oil strata, the following
three approximate formulations are also treated:

\begin{itemize}
\item {\it The lumped formulation.} In this case, stratum's thermal
conductivity is infinitely large in the vertical direction and is
finite in the horizontal direction. The cap and base rock is
considered to be thermally isotropic.

\item {\it The incomplete lumped formulation.} The horizontal heat
transfer in both the cap and base rocks surrounding the strata is
neglected.

\item {\it The Lauwerier formulation.} The horizontal heat transfer
in the strata is also neglected.
\end{itemize}

In the linear case of the incomplete lumped formulation the
temperature field $u = u(x, z, t)$ satisfies the differential
equation given as follows \cite{1}:
\begin{equation}
\frac{\partial u}{\partial t} =
{a}^{2}\frac{\partial^2{u}}{\partial {z}^{2}},\quad 0 < {x, z, t}
< \infty, {a} > {0},
\end{equation}
subject to the boundary condition
\begin{equation}
{z} = {0};\quad \frac{\partial {u}}{\partial {t}} =
\frac{\partial^2{u}}{\partial {x}^{2}} - 2\gamma \frac{\partial
{u}}{\partial {x}} + \alpha \frac{\partial {u}}{\partial
{z}},\quad 0 < {x, t} < \infty,
\end{equation}
and the additional conditions
\begin{align}
 &\hskip -4pc \hbox{(a)}\hskip 2.3pc \quad z = x = 0\hbox{\rm :}\ u = 1,\nonumber\\[.2pc]
 &\hskip -4pc \hbox{(b)}\hskip 2.3pc \quad t = 0\hbox{\rm :}\ u = 0,\nonumber\\[.2pc]
 &\hskip -4pc \hbox{(c)}\hskip 2.3pc \quad u \to 0 \ \hbox{as} \ x^{2} + z^{2}
\to \infty,\phantom{0}\hskip -1pc
\end{align}
where $\gamma$ and $\alpha$ are positive constants.

In the case of the lumped linear formulation, the temperature
field problem is described by (1.1), (1.2) and (1.3) where
(1.3)(a) is replaced by the following condition:
\begin{equation}
{z} = {x} = {0;}\quad \frac{\partial {u}}{\partial {x}} - 2\gamma
{u} = - 2\gamma.
\end{equation}
By using the Laplace and the Hankel transforms, similar boundary
value problems with more general conditions were studied by Ben
Nakhi and Kalla \cite{2}. Some fractional generalizations of the
above discussed problem have been given by Boyadjiev and Scherer
\cite{3,4}. In the present paper we aim at further extending the
works of Ben Nakhi, Kalla, Boyadjiev and Scherer.

In this paper we shall consider the following fractional
differential equation in the linear case of the incomplete lumped
formulation $u = u (x, z, t; \beta)$,
\begin{equation}
{D}_{\ast}^{^{{2}\beta}} {u} = {a}^{2}
\frac{\partial^2{u}}{\partial {z}^{2}} - \lambda \frac{\partial
{u}}{\partial {z}}\hbox{\rm :}\ 0 < {x, z} <  \infty ,\quad 0 <
\beta < \frac{1}{2}
\end{equation}
subject to the condition
\begin{equation}
{z} = {0\hbox{\rm :}\ D}_{\ast}^{{2}\beta} {u} =
\frac{\partial^{2}{u}}{\partial {x}^{2}} - 2\gamma \frac{\partial
{u}}{\partial {x}} + \alpha \frac{\partial {u}}{\partial
{z}},\quad 0 < {x, t} < \infty,
\end{equation}
and the additional conditions
\begin{align}
&\hskip -4pc \hbox{(a)}\hskip 2.3pc\quad z = x = 0\hbox{\rm :}\ u = A h (t),\nonumber\\[.2pc]
&\hskip -4pc \hbox{(b)}\hskip 2.3pc\quad \hbox{if} \ x^{2} + z^{2} \to \infty, \ \hbox{then} \  u \to 0,\nonumber\\[.2pc]
&\hskip -4pc \hbox{(c)}\hskip 2.3pc\quad t = 0\hbox{\rm :}\ u =
0,\phantom{0}\hskip -1pc
\end{align}
where $\lambda, A, \gamma, \alpha$ are positive constants.

We also discuss the linear case of the complete lumped formulation
when the problem remains the same as given by (1.5) subject to the
boundary conditions (1.6) and the conditions (1.7)(b) and (c) with
(1.7)(a) replaced by
\begin{equation}
z = x = 0\hbox{\rm :}\ \frac{\partial {u}}{\partial {x}} - 2\gamma
{u} = - 2\gamma {A}.
\end{equation}
We shall also study the following fractional incomplete lumped
formulation in the radial case. Equation (1.5) in radial
coordinates is considered as subject to the condition
\begin{equation}
{z} = {0\hbox{\rm :}\ D}_{\ast}^{2\beta} {u} =
\frac{\partial^{2}{u}}{\partial {r}^{2}} + \frac{1 - 2\nu}{{r}}
\frac{\partial {u}}{\partial {r}} + \alpha \frac{\partial
{u}}{\partial {z}},\quad 0 < {r, t} <  \infty,
\end{equation}
and the additional conditions
\begin{align}
&\hskip -4pc \hbox{(a)}\hskip 2.3pc\quad z = r = 0\hbox{\rm :}\ u =A h (t),\nonumber\\[.2pc]
&\hskip -4pc \hbox{(b)}\hskip 2.3pc\quad \hbox{if} \ r^{2} + z^{2} \to \infty, \ \hbox{then} \ u \to 0,\nonumber\\[.2pc]
&\hskip -4pc \hbox{(c)}\hskip 2.3pc\quad t = 0\hbox{\rm :}\ u =
0\phantom{0}\hskip -1pc
\end{align}
Here $D_*^{^{2\beta}}$ denotes the Caputo fractional derivative
\cite{9} defined as follows:
\begin{equation}
{D}_{\ast}^{^\beta} {f (t)} =  \left\{ \begin{array}{l}
\dfrac{{1}}{\Gamma {(m} - \beta)}\displaystyle \int_{_0}^{^{t}}
\dfrac{{f}^{m}({x)}}{{(t} - {x)}^{\beta + 1 -
{m}}} {\rm d}x,\quad m - 1 < \beta < {m} \\[1.3pc]
\dfrac{{\rm d}^{m} {f (t)}}{{\rm d}t^{m}},\quad \beta = {m},\quad
{m} \in {N}
\end{array}.\right.
\end{equation}
The following rule of the Laplace transform will be used:
\begin{equation}
{L [D}_{\ast}^{^\beta} {f(t)]} = {p}^\beta {L (f(t)]} -
\sum\limits_{{k} = {0}}^{{m} - 1} {f}^{k}(0) {p}^{\beta - 1 -
{k}},\quad {m} - 1 < \beta \le {m},
\end{equation}
where the Laplace transform of a function $f(t)$ is defined as
\begin{equation}
{\bar {f} (p)} = {L [f(t)]} = \int_{_0}^{^\infty} {\rm e}^{-{pt}}
{f (t) {\rm d}t,\quad {\rm Re} (p)} > {0}.
\end{equation}

\section{Auxiliary results}

\setcounter{equation}{0}

\begin{efros}{\it
If $G(p)$ and $q(p)$ are analytic functions and
\begin{equation}
{F (p)} = L [f(t)], {\rm e}^{-{q(p)}\tau} {G (p)} = {L
[g(t,}\tau)],
\end{equation}
then
\begin{equation}
{F(q(p)) G(p)} = {L}\left[{\int_{_0}^{^\infty} {f (}\tau) {g
(t,}\tau) {\rm d}\tau} \right],
\end{equation}
where $L$ is the Laplace transform operator.}\vspace{.5pc}
\end{efros}

Efros' theorem is a generalization of the Laplace transform
convolution theorem when $q(p) = p$.

The following results will help in solving our problem \cite{3,8}:

If $0 < t, \tau , z < \infty$,
\begin{equation}
{g' (t,}\tau,{z;}\beta) = \frac{({z} + \alpha \tau)\beta}{{a}
{t}^{\beta + 1}} {M}\left({\frac{{z} + \alpha \tau}{{a t}^\beta};
\beta} \right),\quad 0 < \beta < 1
\end{equation}
and
\begin{equation}
{g''(t,}\tau; \beta) = {N}\left({\frac{\tau}{{t}^{2\beta}};\
\beta} \right), \quad 0 < \beta \le \frac{1}{2}
\end{equation}
then
\begin{align}
&\hskip -4pc \hbox{(i)}\hskip 2.3pc\quad {\rm e}^{-\frac{({z} + \alpha \tau )}{{a}} {p}^\beta} = L \{g' (t,\tau, {z;}\beta)\},\\[.4pc]
&\hskip -4pc \hbox{(ii)}\hskip 2.3pc\quad \frac{1}{{p}} {\rm
e}^{-\tau {p}^{{2}\beta}} = L \{g'' (t,\tau;\beta)\}.
\end{align}
The above results follow directly from p.~22, Lemma~2 of \cite{3}
on replacing $b = {\alpha}/{a}$.

Here, functions $M$ and $N$ are auxiliary functions defined as
\begin{equation}
{M (z;}\beta) = \frac{1}{2\pi {i}} \int_{_{H}} \frac{1}{\sigma^{1
- \beta}} {\rm e}^{\sigma - {z}\sigma^\beta} {\rm d}\sigma,\quad 0
< \beta < 1
\end{equation}
and
\begin{equation}
{N (z;}\beta) = \frac{1}{2\pi {i}} \int_{_{H}} \frac{1}{\sigma}
{\rm e}^{\sigma - {z}\sigma^{2\beta}}{\rm d}\sigma,\quad 0 < \beta
< \frac{1}{2},
\end{equation}
where $H$ denotes the Hankel path of integration that begins at
$\sigma = - \infty - ib_{1} (b_{1} > 0)$, encircles the branch
point that lies along the negative real axis and ends up at
$\sigma = - \infty + ib_{2} (b_{2} > 0)$.

The above functions are special cases of the following Wright
function (or Bessel--Maitland function) \cite{5,11}
\begin{equation}
{W (z;}\lambda,\mu) = \sum\limits_{{n} = {0}}^\infty
\frac{{z}^{n}}{{n !} \Gamma (\lambda {n} + \mu)} = \frac{1}{2\pi
{i}} \int_{_{H}} \frac{{\rm e}^{\sigma + {z}
\sigma^{-\lambda}}}{\sigma^\mu}{\rm d}\sigma,\quad \lambda
> - 1, \ \mu > 0
\end{equation}
which is an entire function of $z$.

From eqs~(2.7), (2.8) and (2.9) the following relations can easily
be deduced:
\begin{equation}
{M (z;}\beta) = {W (} - {z;} - \beta, 1 - \beta)
\end{equation}
and
\begin{equation}
{N (z;}\beta) = {W (} - {z;} - 2\beta, 1).
\end{equation}
For $\beta = {1}/{2}$,
\begin{align}
M \left(z;\frac{1}{2}\right) &= \frac{1}{\sqrt \pi} {\rm
e}^{-({z}^{2}/{4})},\\[.4pc]
N \left(z; \frac{{1}}{{2}} \right) &= H (t -z).
\end{align}
Here $H$ denotes the Heaviside step function.

For $\beta = 1$, we have the known result \cite{8}
\begin{equation*}
M(z,1) = \delta (z - 1),
\end{equation*}
where $\delta$ is the Dirac delta function.

\begin{lem}{\it
If $0 < t, \tau, z < \infty,$ then
\begin{equation}
\hskip -4pc \hbox{\rm (a)}\hskip 2.3pc\quad {L}^{-1} ({\rm
e}^{-\tau {p}^{{2}\beta}}) = \frac{2\beta \tau}{{t}^{{2}\beta +
1}} {M}\left({\frac{\tau}{{t}^{{2}\beta}}; 2\beta} \right),\quad 0
< 2\beta < {1}
\end{equation}
and
\begin{align}
&\hskip -4pc \hbox{\rm (b)}\hskip 2.3pc \quad {L}^{-1} \left({\rm
e}^{-\frac{({z} + \alpha \tau)}{{a}} \sqrt {{p}^{{2}\beta} +
\frac{\lambda^2}{{4a}^{2}}}}\right) = \frac{\beta ({z} + \alpha
\tau)}{{a}\sqrt \pi {t}^{{2}\beta + 1}} \int_{_0}^{^\infty}
{u}^{-1 / 2} \cdot {\rm e}^{-\frac{({z} + \alpha
\tau)^2}{{4a}^{2}{u}}}
{\rm e}^{-\frac{\lambda^2{u}}{{4a}^{2}}}\nonumber\\[.6pc]
&\hskip 8.6pc\quad\, \cdot M \left(\frac{{u}}{{t}^{{2}\beta}};
2\beta\right)
{\rm d}u, \quad \hbox{for} \ 0 < {2}\beta \le 1.\\[-3.5pc]\nonumber
\end{align}}
\end{lem}\vspace{.5pc}

\begin{proof} Part (a) can easily be obtained on replacing $\beta$
by $2\beta\ (0 < 2 \beta < 1)$ and $z = 0; a =\alpha = 1$ in (2.3)
and (2.5).

To prove part (b) we apply Efros' theorem. We first express
\begin{equation*}
{\rm e}^{-\frac{({z} + \alpha \tau)}{{a}} \sqrt {{p}^{{2}\beta} +
\frac{\lambda^2}{{4a}^{2}}}} = {F (q(p)) G (p)},
\end{equation*}
where ${G (p)} = {1, q (p)} = {p}^{{2}\beta} +
({\lambda^2}/{{4a}^{2}})$ so that ${F (q(p))} = {\rm
e}^{-\frac{({z} + \alpha \tau)}{{a}} \sqrt {{q(p)}}}$.

Using a known result \cite{6}, we get
\begin{align*}
{F (p)} = {\rm e}^{-\frac{({z} + \alpha \tau )}{{a}} \sqrt {p}} &=
{L}\left\{{\left(\frac{{z} + \alpha \tau }{{a}}\right) \frac{1}{2
\sqrt {\pi {t}^{3}}} {\rm e}^{-\frac{({z} + \alpha
\tau)^2}{{4a}^{2}{t}}}}
\right\}\\[.4pc]
&= {L (f(t))}.
\end{align*}
Hence $G (p) {\rm e}^{-\tau _1 {q(p)}} = {\rm e}^{-\tau _1
({p}^{{2}\beta} + ({\lambda^2}/{{4a}^{2}}))}$. Then
\begin{align*}
{L}^{-1} (G(p) {\rm e}^{-\tau _1 {q(p)}}) &= {\rm
e}^{-\frac{\lambda^2\tau _1}{{4a}^{2}}} \frac{2\beta
\tau_1}{{t}^{{2}\beta + 1}} M \left(\frac{\tau_1}{{t}^{{2}\beta}};
2\beta\right)\\[.4pc]
&= {g (t,}\tau _1).
\end{align*}
By Efros' theorem we can easily get the result (2.15).

In the case of $2\beta = 1$, on using the formula \cite{8}
\begin{equation*}
M (z, 1) = \delta (z -1),
\end{equation*}
where $\delta$ is the Dirac delta function, we get, after some
manipulations and applying the properties of the Dirac delta
function, the following inverse Laplace transform
\begin{equation}
{L}^{-1} \left\{{\rm e}^{-\frac{({z} + \alpha \tau)}{{a}} \sqrt
{{p} + \frac{\lambda^2}{{4a}^{2}}}}\right\} = \frac{({z} + \alpha
\tau)}{{2a}\sqrt \pi} {t}^{-3/2} {\rm e}^{-\frac{({z} + \alpha
\tau)^2}{{4a}^{2}{t}}} {\rm e}^{-\frac{\lambda^2{t}}{{4a}^{2}}}.
\end{equation}
\end{proof}

\section{Fractional extension for incomplete lumped formulation in
the linear case}

\setcounter{equation}{0}

\begin{theor}[\!]
The solution of the linear case of the fractional incomplete
lumped formulation problem $(1.5), (1.6)$ and $(1.7)$ is given by
\begin{equation}
{u (x, z, t)} =  \frac{A{\rm e}^{\gamma {x}}}{2 \sqrt \pi}
\int_{_0}^{^\infty} \frac{x {\rm e}^{-\frac{{x}^{2}}{4\tau} -
\gamma^2\tau}}{\tau^{3 / 2}} {\rm e}^{\frac{\lambda }{{2a}^{2}}
({z} + \alpha \tau )} {g}_{1} ({t,}\tau, {z}; \beta ) {\rm d}\tau,
\end{equation}
where $g_{1}(t, \tau, z; \beta)$ is defined as the Laplace
transform convolution of three functions given by eq.~$(3.17)$.
\end{theor}

\begin{proof}
Taking the Laplace transform of eq.~(1.5), boundary condition
(1.6) and the conditions (1.7) with respect to variable $t$, we
obtain by using the result (1.12),
\begin{align}
&{p}^{{2}\beta} {\bar {u}} = {a}^{2} \frac{\partial^2{\bar
{u}}}{\partial {z}^{2}} - \lambda \frac{\partial {\bar
{u}}}{\partial {z}} \hbox{\rm :}\ 0 < {x, \ z} < \infty,\\[.4pc]
&{z} = {0 \hbox{\rm :}\ p}^{{2}\beta} {\bar {u}} =
\frac{\partial^{2}{\bar {u}}}{\partial {x}^{2}} - 2\gamma
\frac{\partial {\bar {u}}}{\partial {x}} + \alpha \frac{\partial
{\bar {u}}}{\partial {z}}, 0 < {x} <  \infty,\\[.4pc]
&\hskip -4pc \hbox{(a)}\hskip 2.3pc\quad {z} = {x} = {0 \hbox{\rm :}}\ {\bar {u}} = {A \bar {h}} {(p)},\\[.4pc]
&\hskip -4pc \hbox{(b)}\hskip 2.3pc\quad {\bar {u}} \to {0} \
\hbox{if} \ {x}^{2} + {z}^{2} \to \infty,
\end{align}
where ${\bar {u}} = {\bar {u}} {(x, z, p;}\beta) = {L \{u (x,z,t};
\beta); {t} \to p\}$ and $\bar {h} {(p)}$ denotes the Laplace
transform of $h(t)$ which occurs in the boundary condition
(1.7)(a). The solution of (3.2) remains bounded as $z\to \infty$
and is given by
\begin{equation}
{\bar {u}} {(x, z, p;}\beta) = B (x, p) {\rm e}^{-\frac{{z}}{{2}}
{K}},
\end{equation}
where
\begin{equation}
{K} = - \frac{\lambda}{{a}^{2}} + \sqrt {\frac{\lambda^2}{{a}^{4}}
+ \frac{{4p}^{{2}\beta }}{{a}^{2}}}.
\end{equation}
Using (3.6) in eq.~(3.3) to determine function $B (x,p)$ we get
the ordinary differential equation
\begin{equation}
\frac{{\rm d}^{2}{B}}{{\rm d}x^{2}} - 2\gamma \frac{{\rm d}B}{{\rm
d}x} - \left({{p}^{{2}\beta} + \frac{\alpha}{2} {K}} \right) {B} =
{0,}
\end{equation}
under the conditions obtained from (3.4) and (3.5) as
\begin{equation}
{B (0,p)} = {A \bar {h}} {(p)}\quad \hbox{and}\quad
\lim\limits_{{x} \to \infty} {B (x,p)} = 0.
\end{equation}
The solution of (3.8) with (3.9) is given by
\begin{equation}
{B (x,p)} = A \bar {h} (p) {\rm e}^{{x}\left({\gamma - \sqrt
{\gamma ^2 + {p}^{{2}\beta} + \frac{\alpha}{2} {K}}} \right)}.
\end{equation}
Thus the solution given by eq.~(3.6) is expressed as
\begin{equation}
{\bar {u}} {(x, z, p;}\beta) = A {\rm e}^{\gamma {x}} h (p) {\rm
e}^{-\frac{{z}}{{2}} {K}} - {\rm e}^{-{x}\sqrt {\gamma^2 +
{p}^{2\beta} + \frac{\alpha}{2} {K}}}.
\end{equation}
We shall use Efros' theorem to find the inverse Laplace transform
of (3.11) to obtain the solution of our proposed problem.

Express
\begin{equation*}
{\bar {u}} {(x,z,p;}\beta) = {F (q (p;} \beta)) {G (p)},
\end{equation*}
where
\begin{align}
{G (p)} &= A {\rm e}^{\gamma {x}} \bar {h} (p) {\rm
e}^{-\frac{{z}}{{2}} {K}},\nonumber\\[.4pc]
{q (p)} &=  \gamma^2 + {p}^{{2}\beta} + \frac{\alpha}{2} {K},\nonumber\\[.4pc]
{F (p)} &= {\rm e}^{-{x} \sqrt {p}} = {L}\left\{{\frac{x {\rm
e}^{-{x}^{2} / {4t}}}{2 \sqrt {\pi {t}^{3}}}} \right\} = {L
\{f(t)\}},\\[.4pc]
{\rm e}^{-\tau {q (p;} \beta)} {G (p)} &= A {\rm e}^{\gamma {x} -
\gamma^2\tau} {\rm e}^{-\tau {p}^{{2}\beta}} {\bar {h}} {(p)} {\rm
e}^{-\frac{({z} + \alpha \tau )}{2}{K},}\nonumber\\[.4pc]
&= {L (g (t,}\tau,{z;}\beta)).
\end{align}
Then
\begin{equation}
{u (x, z, t;}\beta) = {L}^{-1} ({\bar {u})} = \int_{_0}^{^\infty}
{f (}\tau ) {g (t,}\tau,{z ;}\beta) {\rm d}\tau,
\end{equation}
where functions $f$ and $g$ are defined in (3.12) and (3.13)
respectively.

Now to find $g (t, \tau, z; \beta)$, we write
\begin{equation}
{g (t,}\tau, {z;}\beta) = A {\rm e}^{\gamma {x} - \gamma^2\tau}
{\rm e}^{\frac{\lambda}{{2a}^{2}} ({z} + \alpha \tau)} {g}_{1}
({t,}\tau,{z;}\beta),
\end{equation}
where
\begin{equation}
{g}_{1} ({t,}\tau, {z;} \beta) = {L}^{-1} \left[{{\rm e}^{-\tau
{p}^{{2}\beta}} \bar {h}(p) {\rm e}^{-\frac{({z} + \alpha
\tau)}{{a}} \sqrt {{p}^{{2}\beta} + \frac{\lambda^2}{{4a}^{2}}}}}
\right].
\end{equation}
Using the Laplace transform convolution theorem,
\begin{equation}
{g}_{1} ({t,}\tau,{z;}\beta ) = {L}^{-1} ({\rm e}^{-\tau
{p}^{{2}\beta}}) * {L}^{-1} ({\bar {h}(p))} * {L}^{-1} \left({{\rm
e}^{-(\frac{{z} + \alpha \tau}{{a}}) \sqrt {{p}^{{2}\beta} +
\frac{\lambda^2}{{4a}^{2}}}}} \right).
\end{equation}
Here in (3.17) the functions can be simplified by using results
(2.14), (1.7)(a) and (2.15).

Hence the required solution is given as
\begin{equation}
{u (x,z,t;}\beta) = \frac{A {\rm e}^{\gamma {x}}}{2 \sqrt \pi}
\int_{_0}^{^\infty} \frac{x {\rm e}^{-\frac{{x}^{2}}{{4}\tau} -
\gamma^2\tau}}{\tau^{3 / 2}} {\rm e}^{\frac{\lambda }{{2a}^{2}}
({z} + \alpha \tau )} {g}_{1} ({t,}\tau,{z;}\beta ) {\rm d}\tau.
\end{equation}
\end{proof}

\begin{scase}$\left.\right.$\vspace{.5pc}

\noindent (i) When $2\beta= 1$, Theorem~1 reduces to the problem
studied by Ben Nakhi and Kalla \cite{2} and the solution can be
obtained by representing the function $g_{1}$ in (3.17) as
\begin{align}
{g}_{1} \left({t,}\tau; \frac{1}{2}\right) &= {L}^{-1} ({\rm
e}^{-\tau {p}}) * {L}^{-1} ({\bar {h}(p))} \ast
{L}^{-1}\left({{\rm e}^{-\frac{({z} + \alpha \tau )}{{a}} \sqrt
{{p} + \frac{\lambda^{2}}{{4a}^{2}}}}} \right)\nonumber\\[.3pc]
&= \{\delta ({t} \!-\! \tau) * {h (t)\}} \ast\left\{{\frac{{(z} +
\alpha \tau)}{{2a}\sqrt \pi} {t}^{-3 / 2} {\rm e}^{-\frac{({z} +
\alpha \tau)^2}{{4a}^{2}{t}}}{\rm
e}^{-\frac{\lambda^2{t}}{{4a}^{2}}}} \right\}\nonumber\\[.3pc]
&= \{{h (t} \!-\! \tau){H (t} \!-\! \tau)\} * \left\{{\frac{{(z} +
\alpha \tau)}{{2a}\sqrt \pi} {t}^{-{3 / 2}} {\rm e}^{-\frac{({z} +
\alpha \tau )^2}{{4a}^{2}{t}}} {\rm e}^{-\frac{\lambda^2}{{4a}^{2}}{t}}} \right\}\nonumber\\[.3pc]
&= \frac{{H (t} \!-\! \tau) ({z} + \alpha \tau)}{{2a}\sqrt \pi}
\int_{_0}^{^{{t} - \tau}} {u}^{-{3} / 2} {\rm e}^{-\frac{({z} +
\alpha \tau )^2}{{4a}^{2}{u}}} {\rm
e}^{-\frac{\lambda^2}{{4a}^{2}}{u}}{h(t} \!-\! \tau \!-\! {u) {\rm
d}u}.
\end{align}
On substitution of (3.19) in (3.18), we get
\begin{align}
u \left(x,z,t,\frac{{1}}{{2}}\right) &= \frac{Ax {\rm e}^{\gamma
{x}}}{{4a}\pi} \int_{_0}^{^{t}} \frac{{\rm
e}^{-\frac{{x}^{2}}{4\tau} - \gamma^2\tau}}{\tau^{3 / 2}} {\rm
e}^{\frac{\lambda}{{2a}^{2}} ({z} + \alpha \tau)}\nonumber\\[.5pc]
&\quad\, \cdot \left[{({z} + \alpha \tau) \int_{_0}^{^{{t} -
\tau}} {u}^{-3 /2}{\rm e}^{-\frac{({z} + \alpha \tau
)^2}{{4a}^{2}{u}}} {\rm e}^{-\frac{\lambda^2{u}}{{4a}^{2}}} {h(t}
\!-\! \tau \!-\! {u) {\rm d}u}} \right]{\rm d}\tau.
\end{align}
(ii) On taking $\lambda= 0, A = h (t) = 1$ in (1.5) and (1.7)(a)
Theorem~1 reduces to the problem considered by Boyadjiev and
Scherer \cite{3} and the solution is obtained by representing the
function $g_{1}$ as a convolution of two functions in Theorem~1
given as
\begin{equation}
{g}_{1} ({t,}\tau, {z;}\beta) = {L}^{-1} \left({\frac{{\rm
e}^{-\tau {p}^{{2}\beta}}}{{p}} {\rm e}^{-\frac{({z} + \alpha
\tau)}{{a}}{p}^\beta}} \right).
\end{equation}
From (2.5) and (2.6), we get
\begin{equation}
{g}_{1} ({t,}\tau; {z;}\beta) = {g' (t,}\tau , {z;}\beta) * {g''
(t,}\tau, {z;}\beta).
\end{equation}
If we further take $\beta = 1/2$, it yields the result given in
Antimirov {\it et~al} \cite{1}.
\end{scase}

\section{Fractional lumped formulation in the linear case}

\setcounter{equation}{0}

\begin{theor}[\!]
The solution of the fractional lumped linear formulation $(1.5),
(1.6),$\break $(1.7)${\rm (b),} {\rm (c)} and $(1.8)$ is given by
\begin{align}
{u (x,z,t;}\beta) &= {2}\gamma A {\rm e}^{\gamma {x}}
\int_{_0}^{^\infty} \left({\frac{1}{\sqrt {\pi \tau}} {\rm
e}^{-\big(\frac{{x}^{2}}{4\tau} + \tau \gamma^2\big)} \!-\! \gamma
{\rm e}^{\gamma {x}} {\rm erfc}\left({\frac{{x}}{{2}\sqrt \tau }
\!+\! \gamma \sqrt \tau} \right)} \!\right)\nonumber\\[.5pc]
&\quad\, \cdot {\rm e}^{\frac{\lambda}{{2a}^{2}}({z} + \alpha
\tau)} {g}_{2} ({t,}\tau, {z;}\beta) {\rm d}\tau,
\end{align}
where $g_{2} (t, \tau, z; \beta)$ is defined as the Laplace
transform convolution of two functions given by $(2.6)$ and
$(2.14)$.
\end{theor}

\begin{proof}
Taking the Laplace transform of the equations of the fractional
lumped formulation problem, it reduces to the following boundary
value problem:
\begin{align}
&{p}^{2\beta} {\bar {u}} = {a}^{2} \frac{\partial^2{\bar
{u}}}{\partial {z}^{2}} - \lambda \frac{\partial {\bar
{u}}}{\partial {z}},\quad 0 < {x, z} <  \infty ,\\[.4pc]
&{z} = {0; p}^{{2}\beta}{\bar {u}} =  \frac{\partial^{2}{\bar
{u}}}{\partial {x}^{2}} - 2\gamma \frac{\partial {\bar
{u}}}{\partial {x}} + \alpha \frac{\partial {\bar {u}}}{\partial
{z}};\quad 0 < {x} <  \infty,\\[.4pc]
&\hskip -4pc \hbox{(a)}\hskip 2.3pc\quad {z} = {x} = {0 \hbox{\rm
:}}\ \frac{\partial {\bar {u}}}{\partial {x}} - 2\gamma {\bar {u}} = - \frac{2\gamma}{{p}} {A},\nonumber\\[.4pc]
&\hskip -4pc \hbox{(b)}\hskip 2.3pc\quad \hbox{if} \ {x}^{2} +
{z}^{2} \to \infty \ \hbox{then} \ \bar {u} \to 0.
\end{align}
Proceeding as in Theorem~1, the solution to (4.2) is
\begin{equation}
{\bar {u}} {(x, z, p;}\beta) = B (x, p) {\rm e}^{-\frac{{z}}{{2}}
{K}},
\end{equation}
where
\begin{equation}
{K}=- \frac{\lambda}{{2a}^{2}} + \frac{{2}}{{a}} \sqrt
{{p}^{{2}\beta} + \frac{\lambda ^2}{{4a}^{2}}}.
\end{equation}
It is bounded as $z \to \infty$.

Now at $z = 0$, (4.5) satisfies eq.~(4.3) i.e.
\begin{equation}
\frac{{\rm d}^{2}{B}}{{\rm d}x^{2}} - 2\gamma \frac{{\rm d}B}{{\rm
d}x} - \left[{{p}^{{2}\beta} + \frac{\alpha}{2} {K}} \right] {B} =
{0}.
\end{equation}
As $\lim_{{x} \to \infty} {B (x,p)} = {0}$ and $B$
satisfies (3.8), solution to (4.7) is
\begin{equation}
{B (x,p)} = B (p) {\rm e}^{{x}\left({\gamma - \sqrt {\gamma^2 +
{p}^{{2}\beta} + \frac{\alpha}{2} {K}}} \right)}, {B (p)} = {B
(0,p)}.
\end{equation}
From (4.4)(a) and (4.8),
\begin{equation}
{B (p)} =  \frac{{2} \gamma {A}}{{p}\left({\gamma + \sqrt
{\gamma^2 + {p}^{{2}\beta} + \frac{\alpha}{2}{K}}} \right)}.
\end{equation}
Hence from (4.8),
\begin{equation}
{B (x,p)} =  \frac{{2}\gamma {A}}{{p}\left({\gamma + \sqrt
{\gamma^2 + {p}^{{2}\beta} + \frac{\alpha}{2}{K}}} \right)} {\rm
e}^{{x}\left({\gamma - \sqrt {\gamma^2 + {p}^{{2}\beta} +
\frac{\alpha}{2} {K}}} \right)}
\end{equation}
and from (4.5), ${\bar {u}}$ is expressed as
\begin{equation}
{\bar {u}} {(x,z, p;}\beta) = \frac{2 \gamma {A}}{{p}\left({\gamma
+ \sqrt {\gamma^2 + {p}^{{2}\beta} + \frac{\alpha}{2}{K}}}
\right)} {\rm e}^{{x}\left({\gamma - \sqrt {\gamma^2 +
{p}^{{2}\beta} + \frac{\alpha}{2}{K}}} \right)} {\rm
e}^{-\frac{{z}}{{2}} {K}}.
\end{equation}
Now using Efros' theorem to express ${\bar {u}} {(x,z,p;}\beta) =
{F(q(p;}\beta)){G (p)}$, we take
\begin{align*}
{q (p;}\beta) &= \gamma^2 + {p}^{2\beta} + \frac{\alpha}{2}
{K},\\[.3pc]
{G (p)} &=  \frac{{2}\gamma {A}}{{p}} {e}^{\gamma {x} -
\frac{{z}}{{2}}{K}}
\end{align*}
and
\begin{equation*}
{F (q (p}; \beta)) = \frac{{\rm e}^{-{x}\sqrt {{q(p ;}\beta
{)}}}}{\gamma + \sqrt {{q (p;}\beta)}}.
\end{equation*}
Using formula \cite{6}
\begin{align}
{F (p)} &=  \frac{{\rm e}^{-{x}\sqrt {p} }}{\gamma + \sqrt {p}} =
{L}\left[{\frac{{1}}{\sqrt \pi} {\rm e}^{-\frac{{x}^{2}}{{4t}}} -
\gamma {\rm e}^{\gamma {x} + \gamma ^2{t}} {\rm
erfc}\left({\frac{{x}}{{2}\sqrt {t}} + \gamma \sqrt {t}} \right)}
\right]\nonumber\\[.4pc]
&= {L (f(t)),}
\end{align}
Now
\begin{align}
{\rm e}^{-\tau {q(p;}\beta)} {G (p;}\beta) &= 2 \gamma A {\rm
e}^{\gamma {x} - \tau \gamma ^2} {\rm e}^{\frac{\lambda}{{2a}^{2}}
({z} + \alpha \tau)} \frac{{\rm e}^{-\tau {p}^{2\beta}}}{{p}} {\rm
e}^{-\frac{({z} + \alpha \tau)}{{a}} \sqrt {{p}^{{2}\beta} +
\frac{\lambda^2}{{4a}^{2}}}}\nonumber\\[.5pc]
&= {L [g(t,}\tau,{z}; \beta)]\\[.5pc]
{\rm e}^{-\tau {q(p;}\beta)} {G (p;}\beta) &= {2}\gamma A {\rm
e}^{\gamma {x} - \tau \gamma ^{2}} {\rm
e}^{\frac{\lambda}{{2a}^{2}}({z} + \alpha \tau)} {L [g}_{2}
({t,}\tau,{z}; \beta)],
\end{align}
where
\begin{align}
{g}_{2} ({t,}\tau,{z}; \beta) &= {L}^{-1} \left[{\frac{{\rm
e}^{-\tau {p}^{{2}\beta}}}{{p}} {\rm e}^{-\frac{({z} + \alpha
\tau)}{{a}} \sqrt {{p}^{{2}\beta} + \frac{\lambda^2}{{4a}^{2}}}}}\right]\nonumber\\[.5pc]
&= {L}^{-1} \left({\frac{{\rm e}^{-\tau {p}^{{2}\beta}}}{{p}}}
\right) * {L}^{-1} \left[{{\rm e}^{-\frac{({z} + \alpha
\tau)}{{a}} \sqrt {{p}^{{2}\beta} + \frac{\lambda^2}{{4a}^{2}}}}}
\right].
\end{align}
Functions in the right-hand side of (4.15) are represented by
(2.6) and (2.15) respectively. Thus the solution (4.1) is obtained
by using results (4.12), (4.13), (4.14) and (4.15) in the
following:
\begin{equation*}
{u (x,z,t;}\beta) = \int_{_0}^{^\infty} {f (}\tau) {g}_{2}
{(t,}\tau,{z}; \beta ) {\rm d}\tau.
\end{equation*}
\end{proof}

\begin{scase}$\left.\right.$\vspace{.5pc}

\noindent (i) When $2\beta = 1$, Theorem~2 reduces to the form
studied by Ben Nakhi and Kalla \cite{2}. To derive it we take
$2\beta = 1$ in (4.15) and use the result (2.16) to get
\begin{align*}
{g}_{2} ({t,}\tau,{z;}\beta) &= {L}^{-1} \left({\frac{{\rm
e}^{-\tau {p}}}{{p}}} \right)* {L}^{-1} \left({{\rm
e}^{-\frac{({z} + \alpha \tau )}{{a}}} \sqrt {{p} +
\frac{\lambda^2}{{4a}^{2}}}}\right)
\end{align*}
\begin{align*}
&= {H (t} - \tau) * \left\{{\frac{({z} + \alpha \tau )}{2{a}\sqrt
\pi} {t}^{-{3 / 2}} {\rm e}^{-\frac{({z} + \alpha \tau
)^2}{{4a}^{2}{t}}} {\rm e}^{-\frac{\lambda^2{t}}{{4a}^{2}}}}
\right\}\\[.5pc]
&= {H (t} - \tau) \int_{_0}^{^{{t} - \tau}} \frac{({z} + \alpha
\tau)}{{2a}\sqrt \pi} {u}^{-3 / 2} {\rm e}^{-\frac{({z} + \alpha
\tau)^2}{{4a}^{2}{u}}} {\rm e}^{-\frac{\lambda^2{u}}{{4a}^{2}}}
{\rm d}u.
\end{align*}
Now using the integral given in \cite{1} we get
\begin{align*}
{g}_{2} ({t,}\tau, {z;} \beta) &= \frac{{a} \sqrt \pi}{({z} +
\alpha \tau)} {H (t} \!-\! \tau) \left[{{\rm
e}^{\frac{\lambda}{{2a}^{2}} ({z} + \alpha \tau)} {\rm
erfc}\left({\frac{\lambda \sqrt {{t} \!-\! \tau}}{{2a}} +
\frac{({z} + \alpha \tau)}{{2a}\sqrt {{t} \!-\! \tau}}} \right)} \right.\\[.5pc]
&\quad\, + {\rm e}^{-\frac{\lambda}{{2a}^{2}} ({z} + \alpha \tau)}
{\rm erfc}\left({\frac{{(z} + \alpha \tau )}{{2a} \sqrt {{t} -
\tau}} - \frac{\lambda}{{2a}} \sqrt {{t} - \tau}} \right).
\end{align*}
Substitution of $g_{2} (t, \tau, z; \beta)$ in (4.1) gives us the
desired result.

\noindent (ii) For $\lambda = 0, A = h(t) = 1$ and setting $b =
{\alpha}/{a}$ the problem (1.5), (1.6), (1.7)(b), (c), (1.8) and
the corresponding solution given by (4.1) reduce to the fractional
formulation in the lumped linear case studied in \cite{3}.

Further, in case $\beta = 1/2$, our problem reduces to the form
considered in Antimirov {\it et~al} \cite{1}.
\end{scase}

\section{Fractional incomplete lumped formulation: The radial case}

\setcounter{equation}{0}

\begin{theor}[\!]
Here we consider temperature field problem for fractional
incomplete lumped formulation in the radial case represented by
$(1.5), (1.9)$ and $(1.10)$ and the solution is given by
\begin{equation}
{u (r, z, t;}\beta) = \frac{{A}}{\Gamma (\nu)} \int_{_0}^{^\infty}
{\frac{1}{\tau^{\nu + 1}}} \left({\frac{{r}}{{2}}} \right)^{2\nu}
{\rm e}^{-\frac{{r}^{2}}{4\tau }} {\rm
e}^{\frac{\lambda}{{2a}^{2}} ({z} + \alpha \tau)} {g}_{1}
({t,}\tau,{z;}\beta ){\rm d}\tau,
\end{equation}
where $g_{1} (t, \tau, z; \beta)$ is the Laplace transform
convolution of the three functions given by eq.~$(3.17)$.
\end{theor}

\begin{proof}
Proceeding as in the previous theorem, if we take the Laplace
transform of problem (1.5) and the set of conditions (1.9) and
(1.10), we get the following:
\begin{align}
&{p}^{2\beta} {\bar {u}} = {a}^{2} \frac{\partial ^2{\bar
{u}}}{\partial {z}^{2}} - \lambda \frac{\partial {\bar
{u}}}{\partial {z}},\quad 0 < {r, \ z} < \infty,\\[.5pc]
&{z} = {0} {\hbox{\rm :}} \ {p}^{2\beta} {\bar {u}} =
\frac{\partial^{2}{\bar {u}}}{\partial {r}^{2}} + \frac{1 -
2\nu}{{r}} \frac{\partial {\bar {u}}}{\partial {r}} + \alpha
\frac{\partial {\bar {u}}}{\partial {z}},\quad 0 < {r} <  \infty,\\[.5pc]
&\hskip -4pc \hbox{(a)}\hskip 2.3pc\quad {r} = {z} = {0 \hbox{\rm
:} \ \bar {u}} = {A \bar {h}} {(p)},\nonumber\\[.5pc]
&\hskip -4pc \hbox{(b)}\hskip 2.3pc\quad \hbox{if} \ r^{2} +
{z}^{2} \to \infty \ \hbox{then} \ \bar {u} \to {0}.
\end{align}
The solution of (5.2) which remains bounded as $z\to \infty$ is
\begin{equation}
{\bar {u}} {(r, z, p;}\beta) = B (r,p) {\rm e}^{-\frac{{K}}{{2}}
{z}}, {K} = - \frac{\lambda}{{a}^{2}} + \frac{2}{{a}} \sqrt
{{p}^{{2}\beta} + \frac{\lambda^2}{{4a}^{2}}}.
\end{equation}
To determine the function $B (r, p)$, substitution of (5.5) into
(5.3) gives the following ordinary differential equation:
\begin{equation}
\frac{{\rm d}^{2}{B}}{{\rm d}r^{2}} + \frac{1 - 2\nu}{{r}}
\frac{{\rm d}B}{{\rm d}r} - \left({{p}^{2\beta} + \frac{\alpha}{2}
{K}} \right) {B} = 0
\end{equation}
and the corresponding initial conditions are obtained by using
(5.5) in (5.4) as follows:
\begin{equation}
{B (0,p)} = {A \bar {h}} (p) \quad \hbox{and} \quad
\lim\limits_{{r} \to \infty} {B (r, p)} = 0.
\end{equation}

Now the bounded solution of eq.~(5.6) with the condition (5.7),
using the formula (p.~986, 8.494 (9) of \cite{7}) is given as
\begin{equation}
{B (r, p)} = {B}_{1} ({p)}\left({\frac{\mu {r}}{{2}}} \right)^\nu
{K}_\nu (\mu {r),}
\end{equation}
where $\mu = \sqrt {{p}^{{2}\beta} + ({\alpha}/{2}){K}}$ and
$K_\nu (\mu {r)}$ is the modified Bessel function of the second
kind, of order $\nu$, and $B_{1}(p)$ is a constant which can be
determined with the help of (5.7) and the standard result.
\begin{equation}
{K}_\nu ({z)} \sim \frac{{1}}{{2}} \Gamma (\nu)
\left({\frac{2}{{z}}} \right)^\nu\quad \hbox{as} \ \ z \to {0}
\end{equation}
is given as
\begin{equation}
{B}_{1} ({p)} =  \frac{{2}}{\Gamma (\nu)} {A} {\bar {h} (p)}.
\end{equation}
Now substitution of (5.10) and (5.8) into (5.5) gives the
solution:
\begin{equation}
{\bar {u}} {(r, z, p;}\beta) = \frac{2}{\Gamma (\nu)} {A \bar {h}
(p)}\left({\frac{\mu {r}}{{2}}} \right)^\nu {K}_\nu (\mu r) {\rm
e}^{-\frac{{zK}}{{2}}}.
\end{equation}
To get the solution of our proposed problem, we shall find the
inverse Laplace transform of (5.11) with the help of Efros'
theorem. So, let us first express
\begin{equation}
{\bar {u}} {(r, z, p;}\beta) = {G (p;}\beta) {F (q (p;}\beta)),
\end{equation}
where
\begin{align}
{G (p;}\beta) &= \frac{2}{\Gamma (\nu)} {A \bar {h}} {(p)}
{\rm e}^{-\frac{{z}}{{2}} {K}},\nonumber\\[.5pc]
{q (p}; \beta) &= {p}^{2\beta} + \frac{\alpha}{2} {K}\\[.5pc]
{F (q(p;}\beta)) &= \left({\frac{\mu {r}}{{2}}} \right)^\nu
{K}_\nu (\mu r),\nonumber
\end{align}
and
\begin{equation*}
\mu = \sqrt {{q (p}; \beta)}.
\end{equation*}
Using the well-known result \cite{6}, we get
\begin{align}
{F (p)} &= \left({\frac{{r}\sqrt{p}}{{2}}} \right)^\nu {K}_\nu
({r}\sqrt {p}) = {L} \left[{\frac{{1}}{{2}}
\left({\frac{{r}}{{2}}} \right)^{2\nu} \frac{1}{{t}^{\nu + 1}}
{\rm e}^{-\frac{{r}^{2}}{{4t}}}} \right]\nonumber\\[.5pc]
&= {L [f (t;}\beta]
\end{align}
Now
\begin{align}
{G (p;} \beta) {\rm e}^{-\tau {q (p}; \beta)} &=
\frac{{2A}}{\Gamma (\nu )} \bar {h} (p) {\rm e}^{-\frac{{z}}{{2}}
{K}} {\rm e}^{-\tau \left({{p}^{{2}\beta} + \frac{\alpha}{{2}}
{K}} \right)}\nonumber\\[.5pc]
&= \frac{{2A}}{\Gamma (\nu)} {\rm e}^{\frac{\lambda }{{2a}^{2}}
({z} + \alpha \tau)} {\rm e}^{-\tau {p}^{{2}\beta}}{\bar {h}} (p)
{\rm e}^{-\frac{({z} + \alpha \tau)}{2} \sqrt {{p}^{{2}\beta} +
\frac{\lambda^2}{{4a}^{2}}}}\nonumber\\[.5pc]
&= {L}\left[{\frac{{2A}}{\Gamma (\nu )} {\rm
e}^{\frac{\lambda}{{2a}^{2}} ({z} + \alpha \tau)} {g}_{1}
({t,}\tau,{z}; \beta)} \right]\nonumber\\[.5pc]
&= {L [g (t,}\tau,{z}; \beta].
\end{align}
Then by Efros' theorem, the solution (5.1) is obtained by using
(5.14) and (5.15) into the following
\begin{equation}
{u (r, z, t;}\beta) = \int_{_0}^{^\infty} {f (}\tau; \beta)
{g}_{1} {(t,}\tau,{z}; \beta) {\rm d}\tau.
\end{equation}

Here ${g}_{1}({t,}\tau,{z;}\beta)$ is the Laplace transform
convolution of the three functions given by eq.~(3.17).
\end{proof}

\begin{scase}$\left.\right.$

\noindent (i)
\begin{coro}$\left.\right.$\vspace{.5pc}

\noindent If we take $2\beta = 1$ in Theorem~{\rm 3,} the problem
reduces to the following form{\rm :}
\begin{align}
&\frac{\partial {u}}{\partial {t}} = {a}^{2}
\frac{\partial^2{u}}{\partial {z}^{2}} - \lambda \frac{\partial
{u}}{\partial {z}};\quad 0 < {r, z, t} <  \infty,\\[.5pc]
&{z} = {0 \hbox{\rm :}}\ \frac{\partial {u}}{\partial {t}} =
\frac{\partial ^2{u}}{\partial {r}^{2}} + \frac{1 - 2\nu}{{r}}
\frac{\partial {u}}{\partial {r}} + \alpha \frac{\partial
{u}}{\partial {z}};\quad 0 < {r, t} <  \infty,\\[.5pc]
&\hskip -4pc \hbox{\rm (a)}\hskip 2.3pc \quad r = z = 0 \hbox{\rm :}\ u = A h (t),\nonumber\\[.5pc]
&\hskip -4pc \hbox{\rm (b)}\hskip 2.3pc \quad \hbox{\rm if} \ {r}^{2} + {z}^{2} \to \infty, \ \hbox{\rm then} \ u \to {0},\nonumber\\[.5pc]
&\hskip -4pc \hbox{\rm (c)}\hskip 2.3pc \quad t = 0 \hbox{\rm :}\
u = 0,
\end{align}
and the solution is given by
\begin{equation}
{u (r, z, t;}1 / 2) = \frac{{A}}{\Gamma (\nu)} \int_{_0}^{^\infty}
\left({\frac{{r}^{2}}{{4}\tau}} \right)^\nu \frac{{\rm
e}^{-\frac{{r}^{2}}{4\tau}}}{\tau} {\rm
e}^{\frac{\lambda}{{2a}^{2}} ({z} + \alpha \tau)} {g}_{1}
({t,}\tau,{z;}1 / 2) {\rm d}\tau,
\end{equation}
where
\begin{align*}
{g}_{1} \left({t,}\tau; {z; }\frac{{1}}{{2}}\right) &= {L}^{-1}
\left({{\rm e}^{-\tau {p}} \bar {h}(p) {\rm e}^{-\frac{({z} +
\alpha
\tau)}{{a}} \sqrt {{p} + \frac{\lambda^2}{{4a}^{2}}}}}
\right)\\[.4pc]
&= {L}^{-1} ({\rm e}^{-\tau {p}} {\bar {h}(p)) \ast L}^{-1}
\left({{\rm e}^{-\frac{({z} + \alpha \tau )}{{a}} \sqrt {{p} +
\frac{\lambda^2}{{4a}^{2}}}}} \right).
\end{align*}
\end{coro}
Using
\begin{equation}
{L}^{-1} ({\rm e}^{-\tau {p}} {\bar {h}(p))} = {h (t} - \tau) {H
(t} - \tau)
\end{equation}
and (2.16) the solution is expressed as
\begin{align}
u \left(r,z, t; \frac{{1}}{{2}}\right) &= \frac{{A}}{\Gamma (\nu)}
\int_{_0}^{^{t}} {\frac{1}{\tau}} \left({\frac{{r}^{2}}{4\tau}}
\right)^\nu {\rm e}^{-\frac{{r}^{2}}{4\tau}} {\rm
e}^{\frac{\lambda}{{2a}^{2}} ({z} + \alpha \tau)}\nonumber\\[.5pc]
&\quad\, \left({\int_{_0}^{^{{t} - \tau}} \frac{({z} + \alpha
\tau)}{{2a}\sqrt \pi} {u}^{-{3 / 2}}{h(t} - \tau - u) {\rm
e}^{-\frac{\lambda^2{u}}{{4a}^{2}}} {\rm e}^{-\frac{({z} + \alpha
\tau)^2}{{4a}^{2} {u}}}{{\rm d}u}} \right) {\rm d}\tau.
\end{align}
(ii) On taking $\lambda = 0, A = h(t) = 1$ and setting $b =
{\alpha}/{a}$ in Theorem~3, we obtain the result given by
Boyadjiev and Scherer \cite{3}.

Further for $2\beta = 1$, we get the result given in Antimirov
{\it et~al} \cite{1}.
\end{scase}

\section*{Acknowledgements}

The authors are highly thankful to the anonymous referee for
valuable suggestions which led to the present form of the paper.
The authors are also thankful to Prof.~S~L~Kalla for useful
discussions. The second author (AR) is thankful to the University
Grants Commission, New Delhi for awarding Teacher Research
Fellowship.

\end{document}